\documentstyle[fleqn,12pt]{article}
\begin{document}\pagenumbering{arabic}\setcounter{page}{1}\pagestyle{plain}\baselineskip=18pt

\thispagestyle{empty} \rightline{YTUMB 2006-01, February 2006}
\vspace{1.4cm}

\begin{center}
{\Large\bf Cartan Calculi On The Quantum Superplane}
\end{center}

\vspace{1cm}
\begin{center} Salih \c Celik\footnote{E-mail: sacelik@yildiz.edu.tr}

Yildiz Technical University, Department of Mathematics, \\
34210 DAVUTPASA-Esenler, Istanbul, TURKEY. \end{center}

\vspace{3cm} {\bf Abstract}

Cartan calculi on the extended quantum superplane are given. To
this end, the noncommutative differential calculus on the extended
quantum superplane is extended by introducing inner derivations
and Lie derivatives.

\newpage
\section{Introduction}

Noncommutative geometry [1] has started to play an important role in different fields of mathematical physics over
the past decade. The basic structure giving a direction to the noncommutative geometry is a differential calculus on an
associative algebra. The noncommutative differential geometry of quantum groups was introduced by Woronowicz [2,3]. In this approach the differential calculus on the group is deduced from the properties of the group and it involves functions on the group, differentials, differential forms and derivatives. The other approach, initiated by Wess and Zumino [4], followed Manin's emphasis [5] on the quantum spaces as the primary objects. Differential forms are defined in terms of noncommuting coordinates, and the differential and algebraic properties of quantum groups acting on these spaces are obtained from the properties of the spaces. The natural extension of their scheme to superspace [6] was introduced by Soni in [7] and Chung in [8]. The noncommutative geometry of the quantum superplane was given in [9].

The differential calculus on the quantum superplane similarly involves functions on the superplane, differentials,  differential forms and derivatives.
The exterior derivative is a linear operator {\sf d} acting on $k$-forms and producing $(k+1)$-forms, such that for scalar functions (0-forms) $f$ and $g$ we have
\begin{eqnarray}
{\sf d} (1) & = & 0, \nonumber \\
{\sf d}(f g) & = & ({\sf d} f) g + (-1)^{deg(f)} \, f \, ({\sf d} g)
\end{eqnarray}
where $deg(f) = 0$ for even variables and $deg(f) = 1$ for odd
variables, and for a $k$-form $\omega_1$ and any form $\omega_2$
\begin{eqnarray}
    {\sf d}(\omega_1 \wedge \omega_2) & = & ({\sf d} \omega_1) \wedge \omega_2 + (-1)^k \, \omega_1 \, \wedge ({\sf d} \omega_2).
\end{eqnarray}
A fundamental property of the exterior derivative {\sf d} is
\begin{eqnarray}
    {\sf d} \, {\sf d} = : {\sf d}^2 & = & 0.
\end{eqnarray}

There is a relationship between the exterior derivative and the Lie derivative. To describe the relation between exterior derivative and the Lie derivative, we introduce a new operator: the inner derivation. Hence the differential calculus on the quantum superplane can be extended into a large calculus. We call this new calculus the Cartan calculus. The connection of the inner derivation denoted by
{\bf \textit i}$_a$
and the Lie derivative denoted by ${\cal L}_a$ is given by the Cartan formula:
\begin{eqnarray}
{\cal L}_a & = & {\sf d} \circ {\bf \textit i}_a + {\bf \textit i}_a \circ {\sf d}.
\end{eqnarray}
This and other formulae are explaned in [10-12]. Here we do not give any details. In the related section we shall give a brief overview without much discussion.

The extended calculus on the quantum plane was introduced in [13] using the approach of [10]. In this work we explicitly set up the Cartan calculi on the quantum superplane. In section 2 and 3 we give some information on the Hopf algebra structures of the quantum superplane and its differential calculus which we shall use in order to establish our notions. The differential structures of Type II and III which appeared in subsection 3.1 exist in Ref. 8, but here they are repeated because we need them. In section 4 we present the commutation rules of the inner derivations and the Lie derivatives with functions on the quantum superplane, differentials and partial differentials.

\newpage
\section{Review of Hopf algebra ${\cal A}$}

Elementary properties of the quantum superplane are described in Ref. 9. We state briefly the properties we are going
to need in this work.

\subsection{The algebra of functions on the quantum superplane}

Let us begin with the definition of the coordinate ring of the quantum superplane ${\cal R}^{1|1}_q$. It is well known
that the quantum superplane [6] is defined as an associative algebra generated by two noncommuting coordinates $x$ and $\theta$ with the relations
\begin{eqnarray}
  x \theta & = & q \theta x, \nonumber \\
  \theta^2 & = & 0,
\end{eqnarray}
where $q$ is a nonzero complex deformation parameter. The algebra
of $q$ polynomials will be called the algebra of functions on the
the quantum superplane and will be denoted by ${\cal A}_0 =: A
({\cal R}_q^{1\vert 1})$. In the limit $q \longrightarrow 1$, this
algebra is commutative and can be considered as the algebra of
polynomials ${\cal R}^{1|1}$ over the usual superplane, where $x$
and $\theta$ are the two coordinate functions.

Let ${\cal A}_0 = A ({\cal R}_q^{1|1})$ be a free unital associative algebra generated by even coordinate $x$ and odd
coordinate $\theta$ obeying relations (5). We extend the algebra ${\cal A}_0$ by including the inverse of $x$ which satisfies
\begin{eqnarray}
x \, x^{-1} & = & 1 = x^{-1} \, x \nonumber
\end{eqnarray}
and we denote it by ${\cal A}$. We know that the algebra ${\cal A}$ is a graded Hopf algebra with the following costructures [9].

\subsection{Hopf algebra structure on ${\cal A}$}

The definitions of a coproduct, a counit and a coinverse on the algebra ${\cal A}$ are as follows:

{\bf (1)} the coproduct $\Delta_{\cal A} : {\cal A} \longrightarrow {\cal A} \otimes {\cal A}$ is an algebra homomorphism and is defined by
\begin{eqnarray}
\Delta_{\cal A}(x) & = & x \otimes x, \nonumber \\
\Delta_{\cal A}(\theta) & = & \theta \otimes x + x \otimes \theta.
\end{eqnarray}

{\bf (2)} The counit $\epsilon_{\cal A}$ is an algebra homomorphism from ${\cal A}$ to the complex numbers ${\cal C}$ and is given by
\begin{eqnarray}
\epsilon_{\cal A}(x) & = & 1, \nonumber \\
\epsilon_{\cal A}(\theta) & = & 0.
\end{eqnarray}

{\bf (3)} The antipode $S_{\cal A} : {\cal A} \longrightarrow {\cal A}$ is an algebra antihomomorphism and is given by
\begin{eqnarray}
S_{\cal A}(x) & = & x^{-1}, \nonumber \\
S_{\cal A}(\theta) & = & - x^{-1} \theta x^{-1}.
\end{eqnarray}
These comaps satisfy the Hopf algebra axioms:
\begin{eqnarray}
(\Delta_{\cal A} \otimes \mbox{id}) \circ \Delta_{\cal A} & = &
   (\mbox{id} \otimes \Delta_{\cal A}) \circ \Delta_{\cal A}, \nonumber \\
m_{\cal A} \circ (\epsilon_{\cal A} \otimes \mbox{id}) \circ \Delta_{\cal A}
   & = & m_{\cal A} \circ (\mbox{id} \otimes \epsilon_{\cal A}) \circ \Delta_{\cal A}, \\
m_{\cal A} \circ (S_{\cal A} \otimes \mbox{id}) \circ \Delta_{\cal A} & = &
   \epsilon_{\cal A} = m_{\cal A} \circ (\mbox{id} \otimes S_{\cal A}) \circ \Delta_{\cal A}. \nonumber
\end{eqnarray}
where id denotes the identity map on ${\cal A}$ and $m_{\cal A}$ stands for the algebra product ${\cal A} \otimes {\cal A} \longrightarrow {\cal A}$. The multiplication in ${\cal A} \otimes {\cal A}$ is defined with the rule
\begin{eqnarray}
(A \otimes B) \, (C \otimes D) = (-1)^{deg(B) \, deg(C)} (AC \otimes BD).
\end{eqnarray}

\section{Differential calculi on the quantum superplane}

In this section, we shall build up the noncommutative differential calculi on the quantum superplane with help of the covariance point of view, using the Hopf algebra structure of the quantum superplane [14].

\subsection{Differential algebra}

It is well known that in classical differential calculus, functions commute with differentials. From algebraic point of view, the space of 1-forms is a free bimodule over the algebra of smooth functions generated by the first order differentials and the commutativity shows how its left and right structures are related to each other.

In order to establish a noncommutative differential calculus on the quantum superplane, we assume
that the commutation relations between the coordinates and their differentials are of the following form:
\begin{eqnarray}
x \, {\sf d} x & = & Q \, {\sf d}x \, x, \nonumber \\
x \, {\sf d} \theta & = & Q_{11} \, {\sf d} \theta \, x + Q_{12} \, {\sf d} x \, \theta, \nonumber \\
\theta \, {\sf d} x & = & Q_{21} \, {\sf d} x \, \theta + Q_{22} \, {\sf d} \theta \, x, \\
\theta \, {\sf d} \theta & = & {\sf d} \theta ~\theta. \nonumber
\end{eqnarray}
The coefficients $Q$ and $Q_{ij} \, (1 \le i, j \le 2)$ will be determined using the covariance of the noncommutative differential calculus. We also have assumed that
\begin{eqnarray}
{\sf d} x \wedge {\sf d} \theta & = & Q' \, {\sf d} \theta \wedge {\sf d} x, \qquad {\sf d} x \wedge {\sf d} x = 0,
\end{eqnarray}
where $Q'$ is a parameter that will be described later.

We shall denote the algebra generated with relations (11) by $\Omega^1$ and the algebra generated with relations (12) by $\Omega^2$.

\subsection{Covariance}

We first note that consistency of a differential calculus with commutation relations (5) means that the differential algebra is a graded associative algebra generated by the elements of the set $\{x, \theta, {\sf d} x, {\sf d} \theta\}$. Let $\Omega$ be a free left module over the algebra ${\cal A}$ generated by the elements of this set. So the $\Omega$ has to be generated by $\Omega^0 \cup \Omega^1 \cup \Omega^2$, where $\Omega^0$ is isomorphic to ${\cal A}$. One says that $(\Omega, {\sf d})$ is a {\it first-order differential calculus} over the Hopf algebra $({\cal A}, \Delta_{\cal A},\epsilon_{\cal A}, S_{\cal A})$. We begin with the definitions of a left- and right-covariant bimodule.

{\bf (1)} Let $\Omega$ be a bimodule over ${\cal A}$ and $\Delta^R : \Omega \longrightarrow \Omega \otimes {\cal A}$
be a linear homomorphism. We say that $(\Omega, \Delta^R)$ is a right-covariant bimodule if
\begin{eqnarray}
  \Delta^R(a \rho + \rho' a') & = & \Delta_{\cal A}(a) \Delta^R(\rho) + \Delta^R(\rho') \Delta_{\cal A}(a')
\end{eqnarray}
for all $a, a' \in {\cal A}$ and $\rho, \rho' \in \Omega$, and
\begin{eqnarray}
  (\Delta^R \otimes \mbox{id}) \circ \Delta^R & = &
  (\mbox{id} \otimes \Delta_{\cal A}) \circ \Delta^R, \nonumber \\
  (\mbox{id} \otimes \epsilon) \circ \Delta^R & = & \mbox{id}.
\end{eqnarray}
The action of $\Delta^R$ on the first order differentials is
\begin{eqnarray}
  \Delta^R({\sf d} x) & = & {\sf d}x \otimes x, \nonumber \\
  \Delta^R({\sf d} \theta) & = & {\sf d} \theta \otimes x + {\sf d} x \otimes \theta
\end{eqnarray}
since
\begin{eqnarray}
  \Delta^R({\sf d} a) & = & ({\sf d} \otimes \mbox{id}) \Delta_{\cal A}(a), \qquad
  \forall a \in {\cal A}.
\end{eqnarray}

We now apply the linear map $\Delta^R$ to relations (11):
\begin{eqnarray}
  \Delta^R(x \, {\sf d} x) & = & \Delta_{\cal A}(x) \Delta^R ({\sf d} x) =
    Q \Delta^R ({\sf d} x \, x), \nonumber \\
  \Delta^R(x \, {\sf d} \theta) & = & \Delta^R (Q_{11} \, {\sf d} \theta \, x + Q_{12} \, {\sf d} x \, \theta) +
    (qQ - Q_{11} - qQ_{12}) {\sf d}x \, x \otimes \theta x, \nonumber \\
  \Delta^R(\theta \, {\sf d} x) & = & \Delta^R (Q_{21} \, {\sf d} x \, \theta + Q_{22} \, {\sf d} \theta \, x) -
    (q^{-1} Q + Q_{21} + q^{-1} Q_{22}) {\sf d}x \, x \otimes x \theta, \nonumber \\
  \Delta^R(\theta \, {\sf d} \theta) & = & \Delta^R ({\sf d} \theta \, \theta) +
    (q^{-1} Q_{11} + Q_{22} - 1) {\sf d} \theta x \otimes x \theta  \nonumber \\ & &
    \,\,\, + (Q_{12} + q C_{21} + 1) {\sf d} x \, \theta \otimes \theta x.
\end{eqnarray}
So we must have
\begin{eqnarray}
Q_{11} + q Q_{12} & = & q Q, \qquad Q_{11} + q Q_{22} = q, \nonumber \\
Q_{12} + q Q_{21} & = & - 1, \qquad q Q_{21} + C_{22} = - Q.
\end{eqnarray}

{\bf (2)} Let $\Omega$ be a bimodule over ${\cal A}$ and $\Delta^L : \Omega \longrightarrow {\cal A} \otimes \Omega$
be a linear homomorphism. We say that $(\Omega, \Delta^L)$ is a left-covariant bimodule if
\begin{eqnarray}
  \Delta^L(a \rho + \rho' a') & = & \Delta_{\cal A}(a) \Delta^L(\rho) + \Delta^L(\rho') \Delta_{\cal A}(a')
\end{eqnarray}
for all $a, a' \in {\cal A}$ and $\rho, \rho' \in \Omega$, and
\begin{eqnarray}
  (\Delta_{\cal A} \otimes \mbox{id}) \circ \Delta^L & = & (\mbox{id} \otimes \Delta^L) \circ \Delta^L, \nonumber \\
  (\epsilon \otimes \mbox{id}) \circ \Delta^L & = & \mbox{id}.
\end{eqnarray}
Since
\begin{eqnarray}
  \Delta^L({\sf d} a) & = & (\tau \otimes {\sf d}) \Delta_{\cal A}(a), \qquad \forall a \in {\cal A}
\end{eqnarray}
the action of $\Delta^L$ on the first order differentials gives rise to the relations
\begin{eqnarray}
  \Delta^L({\sf d} x) & = & x \otimes {\sf d} x, \nonumber \\
  \Delta^L({\sf d} \theta) & = & x \otimes {\sf d} \theta - \theta \otimes {\sf d} x.
\end{eqnarray}
Here $\tau : \Omega \longrightarrow \Omega$ is the linear map of degree zero which gives $\tau(a) = (-1)^{deg(a)} \, a$. If we apply $\Delta^L$ to relations (11), we don't have any new relations between $Q$'s. Consequently, we may have three distinct solutions:

{\bf Type I:} includes one deformation parameter
\begin{eqnarray}
    Q_{12} = 0, \quad Q_{22} = 0 \,\, \Longrightarrow \,\,
    \left\{\begin{array}{r@{\quad}ll} Q & = & 1 \\ Q_{11} & = & q \\ Q_{21} & = & - q^{-1} \\ Q' & = & q \end{array} \right.
\end{eqnarray}

{\bf Type II:} includes two parameters
\begin{eqnarray}
    Q_{22} = 0, \quad Q = r \,\, \Longrightarrow \,\,
    \left\{ \begin{array}{r@{\quad}ll} Q_{11} & = & q \\ Q_{12} & = & r - 1 \\ Q_{21} & = & - q^{-1} \, r \\ Q' & = & q \,      r^{-1} \end{array} \right.
\end{eqnarray}

{\bf Type III:} includes two parameters
\begin{eqnarray}
    Q_{12} = 0, \quad Q = p \,\, \Longrightarrow \,\,
    \left\{ \begin{array}{r@{\quad}ll} Q_{11} & = & p \, q \\ Q_{21} & = & - q^{-1} \\ Q_{22} & = & 1 - p \\ Q' & = & p \, q \end{array} \right.
\end{eqnarray}
Hence we see that solution of the form {\bf Type I} is a special case of solution of the form {\bf Type II} for $r=1$ or {\bf Type III} for $p=1$. Therefore it may be omitted. However the {\bf Type I} solution, as we will see in section 4, gives rise to interesting and important results. Therefore it is convenient it as a special type.

Note that it can be checked that the identities (14), (20) and also the following identities are satisfied:
\begin{eqnarray}
  (\mbox{id} \otimes {\sf d}) \Delta_{\cal A}(a) & = & \Delta^L({\sf d} a), \nonumber \\
  ({\sf d} \otimes \mbox{id}) \Delta_{\cal A}(a) & = & \Delta^R({\sf d} a), \\
  (\Delta^L \otimes \mbox{id}) \circ \Delta^R & = & (\mbox{id} \otimes \Delta^R) \circ \Delta^L.\nonumber
\end{eqnarray}

Note that the {\bf Type I} has a special importance. We will see in section 4, that in this case the extended calculus on the quantum superplane has linear commutation relations, for example, the commutation rules of the Lie derivatives with functions do not contain the inner derivations.

We call the ${\cal A}_0 = A ({\cal R}_q^{1|1})$-bimodule generated by ${\sf d} x$, ${\sf d} \theta$ with relations (11) a {\it cotangent bimodule} and denote it by $\Omega(T^*_. {\cal R}_q^{1|1})$. Further work on this and on the tangent bimodule is in progress.

\subsection{Cartan-Maurer one-forms on ${\cal A}$}

In analogy with the left-invariant one-forms on a Lie group in classical differential geometry, one can construct two one-forms using the generators of $\cal A$ as follows [9]:
\begin{eqnarray}
  \omega_x      & = & {\sf d} x \, x^{-1}, \nonumber \\
  \omega_\theta & = & {\sf d} \theta x^{-1} - {\sf d} x \, x^{-1} \theta x^{-1}.
\end{eqnarray}
The commutation relations between the generators of ${\cal A}$ and one-forms are
\begin{eqnarray}
x \, \omega_x & = & Q \, \omega_x \, x, \nonumber \\
x \, \omega_\theta & = & Q_{11} \, \omega_\theta \, x, \nonumber \\
\theta \, \omega_x & = & -  Q \, \omega_x \, \theta + Q_{22} \, \omega_\theta \, x, \\
\theta \, \omega_\theta & = & Q_{11} \, \omega_\theta \, \theta. \nonumber
\end{eqnarray}

The commutation rules of the one-forms $\omega_x$ and $\omega_\theta$ are
\begin{eqnarray}
\omega_x \wedge \omega_\theta & = & \omega_\theta \wedge \omega_x, \nonumber \\
\omega_x \wedge \omega_x & = & 0.
\end{eqnarray}

We denote the algebra of the forms generated by the two elements $\omega_x$ and $\omega_\theta$ by ${\cal W}$. We make the algebra ${\cal W}$ into a graded Hopf algebra with the following co-structures [9]: the coproduct $\Delta_{\cal W} : {\cal W} \longrightarrow {\cal W} \otimes {\cal W}$ is defined by
\begin{eqnarray}
\Delta_{\cal W}(\omega_x) & = & \omega_x \otimes 1 + 1 \otimes \omega_x, \nonumber \\
\Delta_{\cal W}(\omega_\theta) & = & \omega_\theta \otimes 1 + 1 \otimes \omega_\theta.
\end{eqnarray}
The counit $\epsilon_{\cal W} : {\cal W} \longrightarrow {\cal C}$ is given by
\begin{eqnarray}
\epsilon_{\cal W}(\omega_x) & = & 0, \nonumber \\ \epsilon_{\cal W}(\omega_\theta) & = & 0
\end{eqnarray}
and the coinverse $S_{\cal W}: {\cal W} \longrightarrow {\cal W}$ is defined by
\begin{eqnarray}
S_{\cal W}(\omega_x) & = & - \omega_x, \nonumber \\
S_{\cal W}(\omega_\theta) & = & - \omega_\theta.
\end{eqnarray}

\subsection{The algebra of Partial derivatives}

In this section, we introduce commutation relations between the coordinates of the quantum superplane and their partial derivatives. Later, we illustrate the connection between the relations in subsection 3.5.

To proceed, let us obtain the relations of the coordinates with their partial derivatives. We know that the exterior differential {\sf d} can be expressed in the form
\begin{eqnarray}
{\sf d} f & = & ({\sf d} x \, \partial_x + {\sf d} \theta \, \partial_\theta) f.
\end{eqnarray}
Then, for example,
\begin{eqnarray*}
{\sf d} (x f)
& = & {\sf d} x \, f + x \, {\sf d} f\\
& = & [{\sf d} x \, (1 + Q \, x \, \partial_x + Q_{12} \, \theta \, \partial_\theta) + Q_{11} \, {\sf d}            \theta \, \partial_\theta] f \\
& = & ({\sf d} x \, \partial_x \, x + {\sf d} \theta \, \partial_\theta x) f
\end{eqnarray*}
so that, after some calculations
\begin{eqnarray}
\partial_x \, x & = & 1 + Q \, x \, \partial_x + Q_{12} \, \theta \, \partial_\theta, \nonumber \\
\partial_x \, \theta & = & - Q_{21} \, \theta \, \partial_x, \nonumber \\
\partial_\theta \, x & = & Q_{11} \, x \, \partial_\theta,  \\
\partial_\theta \, \theta & = & 1 - \theta \, \partial_\theta - Q_{22} \, x \, \partial_x. \nonumber
\end{eqnarray}
We shall denote the ${\cal A}$-module generated by the partial derivatives $\partial_x$ and $\partial_\theta$ with relations (34) by ${\cal D}^1$. This space ${\cal D}^1$ is the bimodule of first order partial differential operators. 

The commutation relations between derivatives are
\begin{eqnarray}
\partial_x \partial_\theta & = & Q' \, \partial_\theta \partial_x, \nonumber \\
\partial_\theta^2 & = & 0.
\end{eqnarray}

The relations between partial derivatives and differentials are found as
\begin{eqnarray}
\partial_x \, {\sf d} x & = & Q^{-1} \, {\sf d} x \, \partial_x - (1 + Q'^{-1} \, Q_{21}^{-1}) \, {\sf d} \theta \,         \partial_\theta, \nonumber\\
\partial_x \, {\sf d} \theta & = & Q_{11}^{-1} \, {\sf d} \theta \, \partial_x, \nonumber \\
\partial_\theta \, {\sf d} x & = & Q_{21}^{-1} \, {\sf d} x \, \partial_\theta, \\
\partial_\theta \, {\sf d} \theta & = & {\sf d} \theta \, \partial_\theta + (1 - Q' \, Q_{11}^{-1}) \, {\sf d} x \,             \partial_x. \nonumber
\end{eqnarray}

The relations between partial derivatives and the exterior derivative, which guarantee the consistence with the basic requirement for the nilpotency of {\sf d}, are
\begin{eqnarray}
    \partial_x \, {\sf d} & = & Q^{-1} \, {\sf d} \, \partial_x, \nonumber \\
    \partial_\theta \, {\sf d} & = & - Q^{-1} \, {\sf d} \, \partial_\theta.
\end{eqnarray}

\subsection{Quantum Lie superalgebra}

The commutation relations of Cartan-Maurer forms allow us to construct the algebra of the generators. In order to obtain the quantum Lie superalgebra of the algebra generators we first write the Cartan-Maurer forms as
\begin{eqnarray}
{\sf d} x & = & \omega_x \, x, \nonumber \\
{\sf d} \theta & = & \omega_x \theta + \omega_\theta \, x.
\end{eqnarray}
The differential {\sf d} can then the expressed in the form
\begin{eqnarray}
{\sf d} = \omega_x \, H + \omega_\theta \, \nabla.
\end{eqnarray}
Here $H$ and $\nabla$ are the quantum Lie superalgebra generators (vector fields). We now shall obtain the commutation relations of these generators. Considering an arbitrary function $f$ of the coordinates of the quantum superplane and using that ${\sf d}^2 = 0$ and
\begin{eqnarray}
{\sf d} \omega_x & = & 0, \nonumber \\
{\sf d} \omega_\theta & = & 0,
\end{eqnarray}
one has the following commutation relations for the quantum Lie superalgebra:
\begin{eqnarray}
H \, \nabla & = & \nabla \, H, \nonumber \\
\nabla^2 & = & 0.
\end{eqnarray}

The commutation relation (41) of the algebra generators should be consistent with monomials of the coordinates of the quantum superplane. To do this, we evaluate the commutation relations between the generators of algebra and the coordinates. The commuation relations of the generators with the coordinates can be extracted from the graded Leibniz rule:
\begin{eqnarray*}
{\sf d} (x f)
& = & ({\sf d} x) f + x ({\sf d} f) \\
& = & \omega_x \, (x + Q \, x \, H) f + \omega_\theta \, (x \, \nabla) f \\
& = & (\omega_x \, H + \omega_\theta \, \nabla) x f
\end{eqnarray*}
and
\begin{eqnarray*}
{\sf d} (\theta f)
& = & ({\sf d} \theta) f + \theta ({\sf d} f) \\
& = & \omega_x \, (\theta + Q \, \theta \, H) f + \omega_\theta \, (x - Q_{11} \, \theta \, \nabla - Q_{22} \, x \, H) f \\
& = & (\omega_x \, H + \omega_\theta \, \nabla) y f.
\end{eqnarray*}
This yields
\begin{eqnarray}
H \, x  & = & x + Q \, x H, \nonumber \\
H \, \theta & = & \theta + Q \, \theta \, H, \nonumber \\
\nabla \, x & = & Q_{11} \, x \, \nabla, \\
\nabla \, \theta & = & x - Q_{11} \, \theta \, \nabla - Q_{22} \, x \, H. \nonumber
\end{eqnarray}
We know, from subsection 3.4, that the exterior differential {\sf d} can be expressed in the form (33), which we repeat here,
\begin{eqnarray}
{\sf d} f & = & ({\sf d} x \, \partial_x + {\sf d} \theta \, \partial_\theta) f. \nonumber
\end{eqnarray}
Considering (39) together (33) and using (27) one has
\begin{eqnarray}
H & \equiv & x \, \partial_x + \theta \, \partial_\theta, \nonumber \\
\nabla & \equiv & x \, \partial_\theta.
\end{eqnarray}
Note that, using the relations (34) and (35) one can check that the relation of the generators in (43) coincide with (41). It can also be verified that, the action of the generators in (43) on the coordinates coincide with (42).
Of course, these relations can also be found using the dual pairing. This case will be considered in the end of this section.

The commutation relations of the vector fields $H$ and $\nabla$ with the differentials are following
\begin{eqnarray}
    H \, {\sf d} x & = & {\sf d} x \, H, \nonumber \\
    H \, {\sf d} \theta & = & {\sf d} \theta \, H, \nonumber \\
    \nabla \, {\sf d} x & = & Q Q_{21}^{-1} \, {\sf d} x \, \nabla, \\
    \nabla \, {\sf d} \theta & = & Q_{11} \,  {\sf d} \theta \, \nabla + Q_{12} \, {\sf d} x \, H. \nonumber
\end{eqnarray}
Here we used that
\begin{eqnarray}
    Q_{22} - Q_{11} Q_{21} - Q_{11} Q'^{-1} & = & 0, \qquad Q_{12} + Q_{21} (Q_{11} - Q') =  0, \nonumber \\
    Q (Q_{11} - Q') - Q_{11} Q_{12} & = & 0, \\
    Q_{12} (1 + Q' Q_{21}) & = & 0, \qquad Q_{22} (Q_{11} - Q') = 0. \nonumber
\end{eqnarray}

Using relations (44) together with (42) we obtain the commutation rules of the vector fields with one-forms as follows
\begin{eqnarray}
    H \, \omega_x & = & - Q^{-1} \, \omega_x + Q^{-1} \, \omega_x \, H, \nonumber \\
    H \, \omega_\theta & = & - Q^{-1} \, \omega_\theta + Q^{-1} \, \omega_\theta \, H, \nonumber \\
    \nabla \, \omega_x & = & - \omega_x \, \nabla, \\
    \nabla \, \omega_\theta & = & Q^{-1} \, \omega_x + \omega_\theta \, \nabla + (Q - 1) \, \omega_x \, H \nonumber
\end{eqnarray}
or taking
\begin{eqnarray}
    T & = & {\bf 1} + (Q - 1) \, H
\end{eqnarray}
one has
\begin{eqnarray}
    T \, \omega_x & = & Q^{-1} \, \omega_x \, T, \nonumber \\
    T \, \omega_\theta & = & Q^{-1} \, \omega_\theta \, T, \nonumber \\
    \nabla \, \omega_x & = & - \omega_x \, \nabla, \\
    \nabla \, \omega_\theta & = & \omega_\theta \, \nabla + Q^{-1} \, \omega_x \, T.\nonumber
\end{eqnarray}
Here we used that
\begin{eqnarray}
    Q_{12} - Q_{22} = Q - 1.
\end{eqnarray}
Similarly, we can find the commutation relations between the vector fields and the partial derivatives as
\begin{eqnarray}
    \partial_x \, H & = & \partial_x + Q \, H \, \partial_x, \nonumber \\
    \partial_\theta \, H & = & \partial_\theta + Q \, H \, \partial_\theta, \nonumber \\
    \partial_x \, \nabla & = & \partial_\theta + Q Q' \, \nabla \, \partial_x, \\
    \partial_\theta \, \nabla & = & - \nabla \, \partial_\theta.\nonumber
\end{eqnarray}
We here again used that
\begin{eqnarray}
    Q_{12} - Q' Q_{21} & = & 1, \qquad Q_{11} - Q' (Q + Q_{22}) = 0.
\end{eqnarray}

We know that the differential operator {\sf d} satisfies the graded Leibniz rule. Therefore, the generators $H$ and $\nabla$ are endowed with a natural coproduct. To find them, we need to the following commutation relation
\begin{eqnarray}
H x^m & = & \frac{1 - Q^m}{1 - Q} \, x^m + Q^m \, x^m \, H
\end{eqnarray}
where use was made of the first relation of (42). The relation (52) is understood as an operator equation. This implies that when $H$ acts on arbitrary monomials $x^m \theta$,
\begin{eqnarray}
H (x^m \theta) & = & \frac{1 - Q^{m + 1}}{1 - Q} \, (x^m \theta) + Q^{m + 1} \, (x^m \theta) \, H
\end{eqnarray}
from which we obtain
\begin{eqnarray}
H & = & \frac{{\bf 1} - Q^N}{1 - Q}
\end{eqnarray}
where $N$ is a number operator acting on a monomial as
\begin{eqnarray}
N(x^m \theta) & = & (m + 1) x^m \theta.
\end{eqnarray}
We also have
\begin{eqnarray}
\nabla(x^m \theta) & = & Q^m_{11} \, x^{m+ 1}- Q_{11}^{m + 1} \, (x^m \theta) \, \nabla - Q_{11} Q_{22} \, x^{m + 1} \, H.
\end{eqnarray}

So, applying the graded Leibniz rule to the product of functions $f$ and $g$, we write
\begin{eqnarray}
{\sf d} (f g) & = & [(\omega_x \, H + \omega_\theta \, \nabla) f] g + f (\omega_x \, H + \omega_\theta \, \nabla) g
\end{eqnarray}
with help of (39). From the commutation relations of the Cartan-Maurer forms with the coordinates of the quantum superplane, we can compute the corresponding relations of $\omega_x$ and $\omega_\theta$ with functions of the coordinates. From (28) we  have
\begin{eqnarray}
  (x^m \theta) \, \omega_x & = & - Q^{m + 1} \, \omega_x  \, (x^m \theta) + Q^m Q_{22} \, \omega_\theta \, x^{m + 1}, \nonumber\\
  (x^m \theta) \omega_\theta & = & Q_{11}^{m + 1} \, \omega_\theta \, (x^m \theta).
\end{eqnarray}

Inserting (58) in (57) and equating coefficients of the Cartan-Maurer forms, we get, for example,
\begin{eqnarray}
  H(f g) & = & (H f) g + \left({\bf 1} + (Q- 1) \, H\right) \, f \, (H g).
\end{eqnarray}
Consequently, we have the coproducts for {\bf Type II}
\begin{eqnarray}
\Delta(H) & = & H \otimes {\bf 1} + \left({\bf 1} + (r - 1) \, H\right) \otimes H \nonumber \\
\Delta(\nabla) & = & \nabla \otimes {\bf 1} + \left({\bf 1} + (r - 1) \, H\right)^{h_2/h_1} \otimes \nabla,
\end{eqnarray}
where
\begin{eqnarray}
    h_1 & = & \ln r, \qquad h_2 = \ln q,
\end{eqnarray}
or with (54)
\begin{eqnarray}
\Delta(N) & = & N \otimes {\bf 1} + {\bf 1} \otimes N \nonumber \\
\Delta(\nabla) & = & \nabla \otimes {\bf 1} + q^N \otimes \nabla.
\end{eqnarray}
The counit and coinverse may be calculated by using the axioms of Hopf algebra:
\begin{eqnarray}
m(\epsilon \otimes \mbox{id}) \Delta(u) & = & u =  m(\mbox{id} \otimes \epsilon) \Delta(u), \nonumber \\
m(\mbox{id} \otimes S) \Delta(u) & = & \epsilon(u) = m(S \otimes \mbox{id}) \Delta(u).
\end{eqnarray}
So we have
\begin{eqnarray}
\epsilon(N) & = & 0, \nonumber \\
\epsilon(\nabla) & = & 0, \nonumber \\
S(N) & = & - N, \\
S(\nabla) & = & - q^N \, \nabla. \nonumber
\end{eqnarray}

The dual of ${\cal A}$, denoted by ${\cal U}$ is generated by the elements $H$ and $\nabla$ obeying relations (41). Multiplication and comultiplication in ${\cal U}$ can be obtained from the corresponding ones in its dual ${\cal A}$, but in above we get the relevant formulae without using its duality with ${\cal A}$. We now call $<U,a>$ the evaluation of $U$ on $a$ where $U \in {\cal U}$ and $a \in {\cal A}$. Using
\begin{eqnarray}
    <U_1 U_2, a> & = & <U_1 \otimes U_2, \Delta_{\cal A}(a)>, \nonumber\\
    <U, a_1 a_2> & = & <\Delta_{\cal U}(U), a_1 \otimes a_2>
\end{eqnarray}
and
\begin{eqnarray}
    <U, {\bf 1}_{\cal A}> = \epsilon_{\cal U}(U), \qquad
    <{\bf 1}_{\cal U}, a> = \epsilon_{\cal A}(a)
\end{eqnarray}
we then can compute all possible pairings from those between the generators $H$, $\nabla$ and $x$, $\theta$. They are given by for {\bf Type II}
\begin{eqnarray}
    \left\langle T, \left(\matrix{x \cr \theta} \right) \right\rangle = \left(\matrix{r \cr 0} \right), \qquad
    \left\langle \nabla, \left(\matrix{x \cr \theta} \right) \right\rangle = \left(\matrix{0 \cr 1} \right).
\end{eqnarray}
To obtain the (left) action of the elements of ${\cal U}$ on the elements of ${\cal A}$, we now may use the following properties [15]:
\begin{eqnarray}
U[a] & = & a_{(1)} <U, a_{(2)}>, \nonumber \\
U[a f] & = & m \circ \Delta_{\cal U}(U) [a \otimes f]
\end{eqnarray}
where $\Delta_{\cal A}(a) = \sum a_{(1)} \otimes a_{(2)}$ is the coproduct of $a$.
For example, one has
\begin{eqnarray}
    T[x f] & = & m \circ \Delta_{\cal U}(T) [x \otimes f] \nonumber \\
    & = & m (T \otimes T) [x \otimes f] \nonumber \\
    & = & T[x] \, T[f] \nonumber \\
    & = & r \, x \, T[f],
\end{eqnarray}
so that
\begin{eqnarray}
    T \, x & = & r \, x \, T.
\end{eqnarray}
Similarly, we find
\begin{eqnarray}
    T \, \theta & = & r \, \theta \, T, \nonumber \\
    \nabla \, x & = & q \, x \, \nabla, \\
    \nabla \, \theta & = & x - q \theta \, \nabla. \nonumber
\end{eqnarray}

\newpage
\section{Extended calculus on the quantum superplane}

A Lie derivative is a derivation on the algebra of tensor fields over a manifold. The Lie derivative should be defined three ways: on scalar functions, vector fields and tensors.

The Lie derivative can also be defined on differential forms. In this case, it is closely related to the exterior derivative. The exterior derivative and the Lie derivative are set to cover the idea of a derivative in different ways. These differences can be hasped together by introducing the idea of an antiderivation which is called an inner derivation.

\subsection{Inner derivations }

Let us begin with some information about the inner derivations.
Generally, for a smooth vector field $X$ on a manifold the inner derivation, denoted by ${\bf \textit i}_X$, is a linear operator which maps $k$-forms to $(k-1)$-forms. If we define the inner derivation ${\bf \textit i}_X$ on the set of all differential forms on a manifold, we know that ${\bf \textit i}_X$ is an antiderivation of degree $- 1$:
\begin{eqnarray}
    {\bf \textit i}_X (\alpha \wedge \beta) & = & ({\bf \textit i}_X \alpha) \wedge \beta + (-1)^k \, \alpha \wedge ({\bf \textit i}_X \beta)
\end{eqnarray}
where $\alpha$ and $\beta$ are both differential forms. The inner derivation ${\bf \textit i}_X$ acts on 0- and 1-forms as follows:
\begin{eqnarray}
    {\bf \textit i}_X (f) & = & 0, \nonumber \\
    {\bf \textit i}_X ({\sf d} f) & = & X(f).
\end{eqnarray}

We now wish to find the commutation relations between the coordinates $x$, $\theta$ and the inner derivations associated with them. In order to obtain the commutation rules of the coordinates with inner derivations, we shall assume that they are of the following form
\begin{eqnarray}
{\bf \textit i}_x \, x & = & A_1 \, x \, {\bf \textit i}_x + A_2 \, \theta \, {\bf \textit i}_\theta, \nonumber \\
{\bf \textit i}_x \, \theta & = & A_3 \, \theta \, {\bf \textit i}_x + A_4 \, x \, {\bf \textit i}_\theta, \nonumber \\
{\bf \textit i}_\theta \, x & = & A_5 \, x \, {\bf \textit i}_\theta + A_6 \, \theta \, {\bf \textit i}_x, \\
{\bf \textit i}_\theta \, \theta & = & A_7 \, \theta \, {\bf \textit i}_\theta + A_8 \, x \, {\bf \textit i}_x. \nonumber
\end{eqnarray}
The coefficients $A_k$ $(1 \le k \le 8)$ will be determined in terms of the deformation parameters in relations (11). But the use of the relations (5) does not give rise any solution in terms of the parameters $Q$ and $Q_{ij}$ $(1 \le i, j \le 2)$ in (11). Howover, we have, at least, the system
\begin{eqnarray}
    A_4 (A_1 - q A_5) & = & 0, \quad    A_4 (A_3 + q A_7) = 0, \quad    A_2 A_8 = 0, \nonumber\\
    A_8 (A_5 - q A_1) & = & 0, \quad    A_8 (q A_1 + A_7) = 0, \quad    A_4 A_8 = 0.
\end{eqnarray}
To find the coefficients, we need the commutation relations of the inner derivations with the differentials of $x$ and $\theta$. 
Since
\begin{eqnarray}
    {\bf \textit i}_{X_i} ({\sf d} X_j) = \delta_{i j}
\end{eqnarray}
we can assume that the relations between the differentials and the inner derivations are of the following form
\begin{eqnarray}
{\bf \textit i}_x \, {\sf d} x & = & 1 + a_1 \, {\sf d} x \, {\bf \textit i}_x + a_2 \, {\sf d} \theta \, {\bf \textit i}_\theta,\nonumber \\
{\bf \textit i}_x \, {\sf d} \theta & = & a_3 \, {\sf d} \theta \, {\bf \textit i}_x +a_4 \, {\sf d} x \, {\bf \textit i}_\theta,\nonumber \\
{\bf \textit i}_\theta \, {\sf d} x & = & a_5 \, {\sf d} x \, {\bf \textit i}_\theta+a_6 \, {\sf d} \theta \, {\bf \textit i}_x, \\
{\bf \textit i}_\theta \, {\sf d} \theta & = & 1 + a_7 \, {\sf d} \theta \, {\bf \textit i}_\theta + a_8 \, {\sf d} x \, {\bf \textit i}_x. \nonumber
\end{eqnarray}
Using relations (12) we get the system
\begin{eqnarray}
   a_2 (a_5 - Q') & = & 0, \qquad a_5 = Q' \, (1 + a_8), \qquad a_1 = - 1, \nonumber \\
   a_2 (a_7 - Q' a_3) & = & 0, \qquad a_2 = Q' a_3 - 1, \qquad a_2 a_6 = 0.
\end{eqnarray}
Now the use of the relations (11) will give
\begin{eqnarray}
  A_1 & = & Q, \qquad  A_2 = Q_{12}, \qquad A_3 = Q_{21}, \qquad A_4 = 0, \nonumber \\
  A_5 & = & Q_{11}, \qquad  A_6 = 0, \qquad A_7 = 1, \qquad A_8 = Q_{22}
\end{eqnarray}
and some additional relations consisting $A_k$, $a_k$ and
\begin{eqnarray}
    a_6 & = & 0. \nonumber
\end{eqnarray}

To find the remaining parameters $a_k$, this time we assume that the commutation relations of the inner derivations with the partial derivatives $\partial_x$ and $\partial_\theta$ are in the following form
\begin{eqnarray}
{\bf \textit i}_x \, \partial_x & = & B_1 \, \partial_x \, {\bf \textit i}_x + B_2 \, \partial_\theta \, {\bf \textit i}_\theta, \nonumber \\
{\bf \textit i}_x \, \partial_\theta & = & B_3 \, \partial_\theta \, {\bf \textit i}_x + B_4 \, \partial_x \, {\bf \textit i}_\theta, \nonumber \\
{\bf \textit i}_\theta \, \partial_x & = & B_5 \, \partial_x \, {\bf \textit i}_\theta + B_6 \, \partial_\theta \, {\bf \textit i}_x, \\
{\bf \textit i}_\theta \, \partial_\theta & = & B_7 \, \partial_\theta \, {\bf \textit i}_\theta + B_8 \, \partial_x \, {\bf \textit i}_x.
\nonumber \end{eqnarray}
Then using the relations (34) we obtain
\begin{eqnarray}
  B_1 & = & Q^{-1}, \qquad  B_2 = 0, \nonumber \\
  B_3 & = & Q^{-1}_{21}, \qquad B_4 = - Q^{-1}_{1} Q^{-1}_{21} Q_{12}, \nonumber \\
  B_5 & = & Q^{-1}_{11}, \qquad  B_6 = - Q^{-1}_{11} Q^{-1}_{21} Q_{22}, \\
  B_7 & = & 1, \qquad B_8 = 0. \nonumber
\end{eqnarray}

If we demand that the commutation rules of the inner derivations with {\sf d} in the form
\begin{eqnarray}
    {\bf \textit i}_x \circ {\sf d} - F \, {\sf d} \circ {\bf \textit i}_x = \partial_x, \nonumber \\
    {\bf \textit i}_\theta \circ {\sf d} - F' \, {\sf d} \circ {\bf \textit i}_\theta = \partial_\theta,
\end{eqnarray}
one has
\begin{eqnarray}
   F & = & - Q^{-1}, \qquad F' = - F, \nonumber \\
   a_2 & = & - Q^{-1} \, Q_{12}, \qquad a_4 = 0, \qquad a_7 = Q^{-1}, \\
   a_3 & = & - Q^{-1} \, Q_{21}, \qquad a_5 = Q^{-1} \, Q_{11}, \qquad a_8 = Q'^{-1} \, Q^{-1} \, Q_{22}. \nonumber
\end{eqnarray}
Consequently, we have the following commutation relations:
\begin{itemize}

    \item the commutation relations of the inner derivations with $x$ and $\theta$
\begin{eqnarray}
{\bf \textit i}_x \, x & = & Q \, x \, {\bf \textit i}_x + Q_{12} \, \theta \, {\bf \textit i}_\theta, \nonumber \\
{\bf \textit i}_x \, \theta & = & Q_{21} \, \theta \, {\bf \textit i}_x, \nonumber \\
{\bf \textit i}_\theta \, x & = & Q_{11} \, x \, {\bf \textit i}_\theta, \\
{\bf \textit i}_\theta \, \theta & = & \theta \, {\bf \textit i}_\theta + Q_{22} \, x \, {\bf \textit i}_x, \nonumber
\end{eqnarray}

   \item the commutation relations between the differentials and the inner derivations
\begin{eqnarray}
{\bf \textit i}_x \, {\sf d} x &=& 1 - {\sf d} x \, {\bf \textit i}_x- Q^{-1} Q_{12} \, {\sf d} \theta \, {\bf \textit i}_\theta,\nonumber \\
{\bf \textit i}_x \, {\sf d} \theta & = & - Q^{-1} Q_{21} \, {\sf d} \theta \, {\bf \textit i}_x,\nonumber \\
{\bf \textit i}_\theta \, {\sf d} x & = & Q^{-1} Q_{11} \, {\sf d} x \, {\bf \textit i}_\theta,  \\
{\bf \textit i}_\theta \, {\sf d} \theta & = & 1 + Q^{-1} \, {\sf d} \theta \, {\bf \textit i}_\theta + (Q Q')^{-1} Q_{22} \, {\sf d} x \, {\bf \textit i}_x, \nonumber
\end{eqnarray}

\item the relations of the inner derivations with the partial derivatives $\partial_x$ and $\partial_\theta$
\begin{eqnarray}
{\bf \textit i}_x \, \partial_x & = & Q^{-1} \, \partial_x \, {\bf \textit i}_x, \nonumber \\
{\bf \textit i}_x \, \partial_\theta & = & Q_{21}^{-1} \, \partial_\theta \, {\bf \textit i}_x - (Q_{11} Q_{21})^{-1} Q_{12} \,              \partial_x \, {\bf \textit i}_\theta, \nonumber \\
{\bf \textit i}_\theta \, \partial_x & = & Q_{11}^{-1} \, \partial_x \, {\bf \textit i}_\theta - (Q_{11} Q_{21})^{-1} Q_{22} \,              \partial_\theta \, {\bf \textit i}_x, \\
{\bf \textit i}_\theta \, \partial_\theta & = & \partial_\theta \, {\bf \textit i}_\theta. \nonumber
\end{eqnarray}
\end{itemize}

\subsection{Lie derivatives}

We know, from the classical differential geometry, that the Lie derivative ${\cal L}$ can be defined as a linear map from the exterior algebra into itself which takes $k$-forms to $k$-forms. For a 0-form, that is, an ordinary function $f$, the Lie derivative is just the contraction of the exterior derivative with the vector field $X$:
\begin{eqnarray}
    {\cal L}_X f & = & {\bf \textit i}_X \, {\sf d} f.
\end{eqnarray}
For a general differential form, the Lie derivative is likewise a contraction, taking into account the variation in $X$:
\begin{eqnarray}
    {\cal L}_X \, \alpha & = & {\bf \textit i}_X \, {\sf d} \alpha + {\sf d} ({\bf \textit i}_X \alpha).
\end{eqnarray}
The Lie derivative has the following properties. If ${\cal F}(M)$ is the algebra of functions defined on the manifold $M$ then
\begin{eqnarray}
    {\cal L}_X : {\cal F}(M) \longrightarrow {\cal F}(M)
\end{eqnarray}
is a derivation on the algebra ${\cal F}(M)$:
\begin{eqnarray}
    {\cal L}_X (a f + b g) & = & a ({\cal L}_X f) + b ({\cal L}_X g), \nonumber \\
    {\cal L}_X (f g) & = & ({\cal L}_X f) \, g + f \, ({\cal L}_X g),
\end{eqnarray}
where $a$ and $b$ real numbers.

The Lie derivative is a derivation on ${\cal F}(M) \times {\cal V}(M)$ where ${\cal V}(M)$ is the set of vector fields on $M$:
\begin{eqnarray}
    {\cal L}_{X_1} (f X_2) & = & ({\cal L}_{X_1} f) \, X_2 + f \, ({\cal L}_{X_1} X_2).
\end{eqnarray}
The Lie derivative also has an important property when acting on differential forms. If $\alpha$ and $\beta$ are two differential forms on $M$ then
\begin{eqnarray}
    {\cal L}_X (\alpha \wedge \beta) & = & ({\cal L}_X \alpha) \wedge \beta + (-1)^k \, \alpha \wedge ({\cal L}_X \beta)
\end{eqnarray}
where $\alpha$ is a $k$-form.

In this section we now wish to find the commutation rules of the Lie derivatives with functions, {\it i.e.} the elements of the algebra ${\cal A}$, their differentials, etc. For example, the relation of ${\cal L}_x$ with $x$ can be obtained, using relations (74) and (77), as follows:
\begin{eqnarray}
    {\cal L}_x \, x & = & ({\bf \textit i}_x \circ {\sf d} + {\sf d} \circ {\bf \textit i}_x) \, x \nonumber \\
    & = & {\bf \textit i}_x \, {\sf d} x + {\sf d} ({\bf \textit i}_x \, x) \nonumber \\
    & = & 1 + a_1 \, {\sf d} x \, {\bf \textit i}_x + a_2 \, {\sf d} \theta \, {\bf \textit i}_\theta + {\sf d} (A_1 \, x \, {\bf \textit i}_x + A_2 \, \theta \, {\bf \textit i}_\theta) \nonumber\\
        & = & 1 + A_1 \, x \, {\cal L}_x + A_2 \, \theta \, {\cal L}_\theta + (A_1 + a_1) \, {\sf d} x \, {\bf \textit i}_x + (A_2 + a_2) \, {\sf d} \theta \, {\bf \textit i}_\theta \nonumber\\
    & = & 1 + Q \, x \, {\cal L}_x + Q_{12} \, \theta \, {\cal L}_\theta + (Q - 1) \,
      ({\sf d} x \, {\bf \textit i}_x + Q^{-1} Q_{12} \, {\sf d} \theta \, {\bf \textit i}_\theta).
\end{eqnarray}
Similarly, one has
\begin{eqnarray}
{\cal L}_x \, \theta & = & - Q_{21} \, \theta \, {\cal L}_x + Q_{21} (1 - Q^{-1}) \, {\sf d} \theta \, {\bf \textit i}_x,             \nonumber \\
{\cal L}_\theta \, x & = & ({\bf \textit i}_\theta \circ {\sf d} - {\sf d} \circ {\bf \textit i}_\theta) \, x \nonumber\\
& = & Q_{11} \, x \, {\cal L}_\theta + Q_{11} (Q^{-1} - 1) \, {\sf d} x \, {\bf \textit i}_\theta, \\
{\cal L}_\theta \, \theta & = & 1 - \theta \, {\cal L}_\theta - Q_{22} \, x \, {\cal L}_x -
      Q_{22} ((Q Q')^{-1} - 1) \, {\sf d} x \, {\bf \textit i}_x + (Q^{-1} - 1) \, {\sf d} \theta \, {\bf \textit i}_\theta.
\nonumber \end{eqnarray}
The following relations can be obtained from (85)
\begin{eqnarray}
{\cal L}_x \, {\sf d} x & = & {\sf d} x \, {\cal L}_x + Q^{-1} Q_{12} \, {\sf d} \theta \, {\cal L}_\theta, \nonumber\\
{\cal L}_x \, {\sf d} \theta & = & - Q^{-1} Q_{21} \, {\sf d} \theta \, {\cal L}_x, \nonumber\\
{\cal L}_\theta \, {\sf d} x & = & - Q^{-1} Q_{11} \, {\sf d} x \, {\cal L}_\theta, \\
{\cal L}_\theta \, {\sf d} \theta & = & Q^{-1} \, {\sf d} \theta \, {\cal L}_\theta + (Q Q')^{-1} Q_{22} \, {\sf d} x \, {\cal L}_x. \nonumber
\end{eqnarray}

Other commutation relations can be similarly obtained. To complete the description of the above scheme, we get below the remaining commutation relations as follows:
\begin{itemize}
    \item the Lie derivatives and partial derivatives
\begin{eqnarray}
{\cal L}_x \, \partial_x & = & \partial_x \, {\cal L}_x, \nonumber\\
{\cal L}_x \, \partial_\theta & = & - Q Q_{21}^{-1} \, \partial_\theta \, {\cal L}_x + Q (Q_{11} Q_{21})^{-1} Q_{12} \partial_x \, {\cal L}_\theta, \nonumber\\
{\cal L}_\theta \, \partial_x & = & Q Q_{11}^{-1} \, \partial_x \, {\cal L}_\theta - Q (Q_{11} Q_{21})^{-1} Q_{22} \partial_\theta \, {\cal L}_x, \\
{\cal L}_\theta \, \partial_\theta & = & - Q \, \partial_\theta \, {\cal L}_\theta.\nonumber
\end{eqnarray}
  \item the inner derivations
\begin{eqnarray}
{\bf \textit i}_x \, {\bf \textit i}_\theta & = & - Q_{11} (Q_{12} - Q)^{-1} \, {\bf \textit i}_\theta \, {\bf \textit i}_x, \nonumber\\
{\bf \textit i}_x \, {\bf \textit i}_x & = & 0.
\end{eqnarray}
  \item the Lie derivatives and the inner derivations
\begin{eqnarray}
{\cal L}_x \, {\bf \textit i}_x & = & {\bf \textit i}_x \, {\cal L}_x, \nonumber\\
{\cal L}_x \, {\bf \textit i}_\theta & = & - Q Q_{21}^{-1} \, {\bf \textit i}_\theta \, {\cal L}_x + Q_{12} (Q - Q_{12})^{-1} \,        {\bf \textit i}_x \, {\cal L}_\theta, \nonumber\\
{\cal L}_\theta \, {\bf \textit i}_x & = & - Q Q_{11}^{-1} \, {\bf \textit i}_x \, {\cal L}_\theta - (Q' Q_{21})^{-1} Q Q_{22} \, {\bf \textit i}_\theta \, {\cal L}_x, \nonumber \\
{\cal L}_\theta \, {\bf \textit i}_\theta & = & Q^{-1} \, {\bf \textit i}_\theta \, {\cal L}_\theta.
 \end{eqnarray}
  \item the Lie derivatives
\begin{eqnarray}
{\cal L}_x \, {\cal L}_\theta & = & Q_{21}^{-1} (Q_{12} - Q) \, {\cal L}_\theta \, {\cal L}_x, \nonumber \\
{\cal L}_\theta^2 & = & 0.
\end{eqnarray}
\end{itemize}

Note that the Lie derivatives can be written as follows:
\begin{eqnarray}
    {\cal L}_x & = & \partial_x + (1 - Q^{-1}) \, {\sf d} \circ {\bf \textit i}_x, \nonumber \\
    {\cal L}_\theta & = & \partial_\theta - (1 - Q^{-1}) \, {\sf d} \circ {\bf \textit i}_\theta,
\end{eqnarray}
or in terms of vector fields and coordinates
\begin{eqnarray}
    {\cal L}_x & = & x^{-1} \, H - x^{-1} \theta x^{-1} \, \nabla + (1 - Q^{-1}) \, {\sf d} \circ {\bf \textit i}_x, \nonumber \\
    {\cal L}_\theta & = & x^{-1} \, \nabla - (1 - Q^{-1}) \, {\sf d} \circ {\bf \textit i}_\theta.
\end{eqnarray}

\noindent
{\bf Acknowledgment}

\noindent
This work was supported in part by {\bf TBTAK} the Turkish Scientific and Technical Research Council.


\end{document}